\documentclass[10pt]{amsart}

	\usepackage{amsmath,amsthm}
	\usepackage{amssymb}
	\usepackage{verbatim}
	\usepackage{enumerate}
	\usepackage{amscd} 
	\usepackage[all]{xy}

	\DeclareMathAlphabet{\mathpzc}{OT1}{pzc}{m}{it}

	\usepackage{setspace}	
	\usepackage{indentfirst} 
	\usepackage{amsthm}	
	\usepackage{graphicx}	
	\usepackage{amssymb}
	\usepackage{mathrsfs}
	
	\numberwithin{equation}{section}
	\newcommand{\Number}{section} 	
	\theoremstyle{plain}
		\newtheorem{thm}{Theorem}[\Number]
		\newtheorem{lem}[thm]{Lemma}
		
		\newtheorem{prp}[thm]{Proposition}
		\newtheorem{cor}[thm]{Corollary}

		\newtheorem{ntn}[thm]{Notation}

		\newtheorem{dfn}[thm]{Definition}
\newtheorem*{UnNumThm}{Theorem}

\newcommand{\beq}{\begin{equation}} 
\newcommand{\eeq}{\end{equation}} 
\newcommand{\beqr}{\begin{eqnarray*}} 
\newcommand{\eeqr}{\end{eqnarray*}} 
\newcommand{\bal}{\begin{align*}} 
\newcommand{\eal}{\end{align*}} 
\newcommand{\bei}{\begin{itemize}} 
\newcommand{\eei}{\end{itemize}} 
 
\newcommand{\alignInd}{\hspace{.25 in} }

\newcommand{\af}{\alpha}

\newcommand{\dt}{\delta} 
\newcommand{\ep}{\varepsilon}

\newcommand{\io}{\iota}

\newcommand{\sm}{\sigma} 
 
\newcommand{\ph}{\varphi} 
\newcommand{\ps}{\psi}

\newcommand{\ta}{\tau}

\newcommand{\Gm}{\Gamma}

\newcommand{\Z}{{\mathbb{Z}}}

\newcommand{\N}{{\mathbb{N}}}

\newcommand{\id}{{\mathrm{id}}}

\newcommand{\Aut}{{\mathrm{Aut}}}

\newcommand{\Mi}{M_{\infty}}

\newcommand{\andeqn}{\,\,\,\,\,\, {\mbox{and}} \,\,\,\,\,\,} 
 
\newcommand{\ts}[1]{{\textstyle{#1}}}

\newcommand{\ca}{$C^*$-algebra} 
 
\newcommand{\pj}{projection}

\newcommand{\hm}{homomorphism}

\newcommand{\ifo}{if and only if}

\newcommand{\mops}{mutually orthogonal \pj s}

\newcommand{\hsa}{hereditary subalgebra} 
 
\newcommand{\mvnt}{Murray-von Neumann equivalent} 
 
\newcommand{\tRp}{tracial Rokhlin property}

\newcommand{\idsfsuca}{infinite dimensional 
stably finite simple unital \ca} 
\newcommand{\idsuca}{infinite dimensional simple unital \ca}

\renewcommand{\S}{\subset}

\title[Crossed products and the tracial Rokhlin property]{Crossed product C*-algebras by finite group actions with the tracial Rokhlin property } 

\author{Dawn Archey}

\date{\today}   

\address{Department of Mathematics\\
Ben Gurion University of the Negev\\
P.O.B 653 Beer Sheva 84105, ISRAEL}
\email[]{archey@math.bgu.ac.il}


\thanks{The author is indebted to the Fields Institute for Research in Mathematical Sciences for its hospitality during the Thematic Program on Operator Algebras in the fall of 2007 and to the Center for Advanced Studies in Mathematics at Ben Gurion University.  This research is part of the author's Ph.D. thesis at the University of Oregon,
completed under the direction of N. Christopher Phillips.  The author was partially supported by NSF grant DMS 0302401.}

\begin{document}
	
\begin{abstract}
Let
$A$ be a stably finite simple unital $C^*$-algebra and suppose $\alpha $ is an action of a finite group $G$ with the tracial Rokhlin property.  Suppose further $A$ has real rank zero and the order on projections over $A$ is determined by traces.  Then the crossed product $C^*$-algebra $C^*(G,A, \alpha)$ also has real rank zero and order on projections over $A$ is determined by traces.  Moreover, if $A$ also has stable rank one, then $C^*(G,A, \alpha)$ also has stable rank one.

\end{abstract}

\maketitle




	\setlength{\parindent}{.5in}

	\pagenumbering{arabic}	
	\setcounter{page}{1}		

\section{INTRODUCTION}
The \tRp \ for finite group actions was introduced in \cite{Ph31} to address the fact that some \ca s can have no actions with the Rokhlin property of certain groups due to K-theoretic obstructions.  The \tRp \ was further developed in \cite{Ph30} .  In \cite{Ph31} and \cite{ELPW}, examples of actions with the \tRp \ are studied.   This paper studies permanence properties of crossed products by actions
with the tracial Rokhlin property.  The goal of this paper  is to prove  Theorem \ref{K0OrdDetByTraces}, Theorem \ref{chp3:thm1} and Theorem \ref{chp4:thm1} which can be collectively summarized in the following theorem.

\begin{UnNumThm}
Let $A$ be an infinite dimensional stably finite simple unital $C^*$-algebra with real rank zero, 
and suppose that the order on projections over $A$ is determined by traces.  Let 
$\alpha \colon G \rightarrow \Aut (A)$ be an action of a finite group with the tracial Rokhlin 
property.  Then the order on projections over $C^*(G, A, \alpha)$ is determined by traces 
and $C^*(G,A, \alpha )$ has real rank zero.  Moreover, if $A$ also has stable rank one, then $C^*(G, A, \alpha)$ has stable rank one.
\end{UnNumThm}

These three theorems are finite group analogs of known results about actions of $\Z$; see \cite{OP1}.  The proof techniques are similar to those used there.  The most significant changes are in Lemma \ref{Basic1} (which is analogous to Lemma 2.5 of \cite{OP1}). This is the key lemma in the proof of all three theorems.  In that proof a different construction of the isomorphism was needed.

It is known that some hypotheses are needed to ensure that the properties of $A$ pass to the crossed product.  
For the real rank zero situation, consider Example 9 in \cite{E1993} of an outer
action of $\Z / 2\Z $ on a simple nuclear C*-algebra with real rank
zero such that the crossed product does not have real rank zero. 

However, partial results are known, such as Theorem 2.6 of \cite{Ph30}.  This theorem says that if a finite group acts on an infinite dimensional simple unital \ca \ with tracial rank zero via an action with the \tRp, then the crossed product has tracial rank zero.

The situation for stable rank one is similar. Various partial results are known such as theorem 4.6 of \cite{JOPT} which gives cancellation of projections (a condition strictly weaker than stable rank one) under fairly mild hypotheses.   Another result along these same lines is Corollary 5.6 of \cite{JO}.  There, the group must be a finite semidirect product of finite abelian groups, but the action is unrestricted.  In that case, if the group acts on a simple unital \ca \ with property (SP) and stable rank one in such a way that the crossed product has real rank zero, then the crossed product has stable rank one.  

As in the real rank zero situation, it is also known that some conditions will be need to ensure stable rank one in the crossed product.  There is a (nonsimple) unital \ca \ $A$ of stable rank one and an action $\Z / 2\Z$ on $A$ such that the crossed product does not have stable rank one.  See Example 8.2.1 of \cite{B14}.  However, there is no known example for a simple \ca .

This paper is organized as follows.  In Section \ref{Sec:tRp} we establish notation, give definitions, and prove basic facts about the \tRp .  In Section \ref{Sec:Key} we prove the key lemmas needed for the proofs of the three main theorems, which are proved in Sections \ref{Sec:OrdPj} (order on projections is determined by traces), \ref{Sec:RR} (real rank zero), and \ref{Sec:tsr} (stable rank one).





\section{THE TRACIAL ROKHLIN PROPERTY}\label{Sec:tRp}

\begin{ntn}
For any projections $p$ and $q$ in $A$, we write $ p \sim q $ if $p$ is {\bf{(Murray-von Neumann) equivalent}} to $q$.   We write $p \precsim q$ if $p$ is {\bf{(Murray-von Neumann) subequivalent}} to $q$.  
\end{ntn}

\begin{ntn}\label{TraceNtn}
Let $A$ be a unital \ca.
We denote by $T (A)$ the set of all tracial states on $A,$
equipped with the weak* topology.
For any element of $T (A),$
we use the same letter for its standard extension to $M_n (A)$
for arbitrary $n,$
and to $\Mi (A) = \bigcup_{n = 1}^{\infty} M_n (A)$ (no closure).
\end{ntn}

\begin{dfn}
Let $A$ be a unital \ca.
We say that the {\emph{order on \pj s over $A$ is determined by traces}}
if whenever $p, q \in \Mi (A)$ are \pj s such that
$\ta (p) < \ta (q)$ for all $\ta \in T (A),$
then $p \precsim q.$
\end{dfn}

\begin{dfn}\label{NewTRPDfn}
Let $A$ be an \idsuca,
and let $\af \colon G \to \Aut (A)$
be an action of a finite group $G$ on $A.$
We say that $\af$ has the
{\emph{tracial Rokhlin property}} if for every finite set
$F \S A,$ every $\ep > 0,$
and every positive element $x \in A$ with $\| x \| = 1,$
there are \mops\  $e_g \in A$ for $g \in G$ such that:
\begin{enumerate}
\item
$\| \af_g (e_h) - e_{g h} \| < \ep$ for all $g, h \in G.$
\item
$\| e_g a - a e_g \| < \ep$ for all $g \in G$ and all $a \in F.$
\item
With $e = \sum_{g \in G} e_g,$ the \pj\  $1 - e$ is \mvnt\  to a
\pj\  in the \hsa\  of $A$ generated by $x.$
\item
With $e$ as in~(3), we have $\| e x e \| > 1 - \ep.$
\end{enumerate}
\end{dfn}

When $A$ is finite, as was shown in Lemma 1.12 of \cite{Ph30},
Condition~(4) of Definition~\ref{NewTRPDfn} is not needed.

The following lemma is the finite group analog of Lemma 1.4 in \cite{OP1} and can be proved using the same argument.

\begin{lem}\label{DfnUsingTrace}
Let $A$ be an \idsfsuca\  with real rank zero and such that the order on \pj s over $A$ is determined by traces.  
Suppose $\af : G \rightarrow A$ is an action of a finite group on $A$.
Then $\af$ has the \tRp\  %
\ifo\  for every finite set $F \S A$ and every $\ep > 0$
there are \mops\  $e_g\in A$ for each $g \in G$ such that:
\begin{enumerate}
\item
$\| \af_h (e_g) - e_{gh} \| < \ep$ for $g\in G.$
\item
$\| e_g a - a e_g \| < \ep$ for all $g \in G$ and all $a \in F.$
\item
With $e = \sum_{g \in G} e_g,$ we have $\ta (1 - e) < \ep$
for all $\ta \in T (A).$
\end{enumerate}
\end{lem}

The next lemma is Lemma 2.4 of \cite{OP1}.

\begin{prp}\label{OrdPjCrPrd} 
Let $A$ be a simple unital infinite dimensional \ca\ %
with real rank zero, 
and assume that the order on \pj s over $A$ is determined by traces. 
Let $\af \colon \Gm \to \Aut (A)$ be an action of a countable 
amenable group. 
Let $p, \, q \in \Mi (A)$ be \pj s such that 
$\ta (p) < \ta (q)$ for every $\Gm$-invariant tracial 
state $\ta$ on $A.$ 
(We extend $\ta$ to $\Mi (A)$ as in Notation \ref{TraceNtn}).
Then there is $s \in \Mi ( C^* (\Gm, A, \af))$ 
such that 
\[ 
s^* s = p, \,\,\,\,\,\, s s^* \leq q, \andeqn s s^* \in \Mi (A). 
\] 
In particular, $p \precsim q$ in $\Mi ( C^* (\Gm, A, \af)).$ 
\end{prp}


\section{THE PIVOTAL LEMMAS}\label{Sec:Key}

The next two lemmas are used as the key lemmas in the proofs of Theorem \ref{K0OrdDetByTraces}, Theorem \ref{chp3:thm1}, and Theorem \ref{chp4:thm1}.  Lemma \ref{Basic1} is the finite group analog of Lemma 2.5 in \cite{OP1}.

\begin{lem}\label{Basic1} 
Let $A$ be an \idsfsuca\ with real rank zero such that the order on 
\pj s over $A$ is determined by traces. Let $G$ be a finite group of order $n$ and let $\af : G\rightarrow  \Aut (A)$ be an action of $G$ with the \tRp. 
Let $\io \colon A \to C^* (G, A, \af)$ be the inclusion map. 
Then for every finite set $F \S C^* (G, A, \af),$ 
every $\ep > 0,$ 
every $N \in \N,$ and
every nonzero positive element $z \in C^* (G, A, \af),$ 
there exist a \pj\ $e \in A \S C^* (G, A, \af),$ 
a unital subalgebra $D \S e C^* (G, A, \af) e,$ 
a \pj\ $f \in A,$ and an isomorphism 
$\ph \colon M_n \otimes f A f \to D,$ such that: 
\begin{enumerate} 
\item 
With $( e_{g, h} )$ for $g,h \in G$ being a system of matrix units for $M_n,$ 
we have $\ph (e_{1, 1} \otimes a) = \io (a)$ for all $a \in f A f$ 
and $\ph (e_{g, g} \otimes 1) \in \io (A)$ for $g\in G.$ 

\item 
With $( e_{g, g} )$ as in~(1), we have 
$\| \ph (e_{g, g} \otimes a) - \io (\af_{g} (a)) \| 
\leq \ep \| a \|$ 
for all $a \in f A f.$ 

\item 
For every $a \in F$ there exist $b_1, \, b_2 \in D$ such that 
$\| e a - b_1 \| < \ep,$ $\| a e - b_2 \| < \ep,$ 
and $\| b_1 \|, \, \| b_2 \| \leq \| a \|.$ 

\item
$e = \sum_{g\in G} \ph (e_{g, g} \otimes 1).$ 

\item
The \pj\ $1 - e$ is \mvnt\ in $C^* (G, A, \af)$ to a 
\pj\ in the \hsa\ of $C^* (G, A, \af)$ generated by $z.$ 

\item 
There are $N$ 
\mops\ $f_1, f_2, \ldots, f_N \in e D e,$ 
each of which is \mvnt\ in $C^* (G, A, \af)$ to $1 - e.$ 
\end{enumerate} 
\end{lem}

\begin{proof} 
We first note, using the same argument as in the proof of Lemma 2.5 of \cite{OP1}, that it is not necessary to check the estimates 
$\| b_1 \|, \, \| b_2 \| \leq \| a \|$ 
in Condition~(3) of the conclusion. 

Now we do the main part of the proof. 
Let $\ep > 0,$ and let $F \S C^* (G, A, \af)$ be a finite set. 
Let $N \in \N,$ and let $z \in C^* (G, A, \af)$ be a nonzero 
positive element.   

Let $u_g$ for $g \in G$ be the standard unitaries  
in the crossed product $C^* (G, A, \af)$.  
We regard $A$ as a subalgebra of $C^* (G, A, \af)$ in the usual way. 

For each $x \in F$ write $x=\sum_{g\in G} a_g u_g .$
Let $S \S A$ be a finite set which contains all the 
coefficients used for all elements of $F.$ 
Let $M = 1 + \sup_{a \in S} \| a \|.$

Let $\delta_0 < \frac{\ep }{16n^2 M}$.  
Let $\delta_1$ be such that if $p_1, p_2$ are projections in a C*-algebra $B$ 
and if $a \in B$ is such that $\|a^* a -p_1 \| \leq \delta_1$ and $\|a a^*-p_2\|  \leq \delta_1$, 
then there is a partial isometry $s \in B$ such that 
$s^* s= p_1$, $s s^* = p_2$, and $\|a-s\| \leq \delta_0.$  
Let $0 < \delta<\min \{ \delta_0, \delta_1, \frac{\ep}{6 n^3}, 1\}$.

Since $A$ has real rank zero
and (by Lemma 1.5 of \cite{Ph30}) 
$\af_g$ is outer for all $g \in G,$ 
Theorem~4.2 of \cite{JO} (with $N = \{ 1 \}$) 
supplies a nonzero \pj\ $q \in A$ which is \mvnt\ in 
$C^* (G, A, \af)$ to a \pj\ in 
${\overline{z C^* (G, A, \af) z}}.$ 
Moreover, Lemma 2.3 of \cite{OP1} provides 
nonzero orthogonal \mvnt\ \pj s 
$q_0, \, q_1, \, \ldots, \, q_{2 N} \in q A q.$

Apply the \tRp\ (Definition~\ref{NewTRPDfn}) with 
$\delta$ in place of $\ep,$ 
with $S$ in place of $F,$ and with $q_0$ in place of $x.$ 
Call the resulting \pj s $e_g$ for each $g \in G$, 
and let $e = \sum_{g \in G} e_g.$ 

Set $f = e_1,$ and define $w_{g,h}=u_{gh^ {-1}}e_h$.  
We claim that the elements $(w_{g,h})_{g,h \in G}$ form a 
$\delta$-approximate system of $n\times n$ matrix units. 
To prove the claim we compute:
\[ \|w_{g,h}^*-w_{h,g}\| =\|u_{gh^{-1}}  e_h u_{gh^{-1}}^*-e_g\| =\| \alpha_{g h^{-1}} (e_h) -e_g\| < \delta.\]
Then, using $e_g e_h =\delta_{g,h} e_h$ we find
\[ 
\|w_{g_1, h_1}  w_{g_2, h_2} - \delta_{g_2, h_1} w_{g_1, h_2}\| 
< \delta . 
\]
For the final condition, since $\|e q_0 e \| > 1- \delta >0$, 
the projection $e$ is nonzero, so $e_g$ is nonzero for each $g \in G$.   
In particular $\|w_{1,1}\|=\|e_1\|=1 > 1- \delta .$  This proves the claim.

Since $(w_{g,h})_{g,h \in G }$ forms a $\delta$-approximate system of matrix units, 
$w_{g,1}$ is an approximate partial isometry for each $g \in G$.  
More specifically, 

\[\| w_{g,1} w_{g,1}^*-e_g\|=\|u_g e_1 e_1 u_g^* -e_g\|=\| \alpha_g(e_1)-e_g\| < \delta\]
since the $e_g$
are the tracial Rokhlin projections.  Also, 
\[ \| w_{g,1}^* w_{g,1}-e_1\|=\|e_1u_g^*u_g e_1-e_1\|=\|e_1-e_1\|=0 < \delta \] 
because $u_g$ is a unitary for each $g$.

Since $\delta < \delta_1$, by the choice of $\delta_1$ 
there exist partial isometries $z_g \in C^*(G,A,\alpha)$ for each $g \in G$ such that 
$\|z_g-w_{g,1}\|< \delta_0$ and such that $z_g z_g^*=e_g$ and $z_g^* z_g=e_1$.  
Moreover, one may check that we may take $z_1=e_1$.

Let $(e_{g,h})_{g,h \in G}$ be an $n \times n$ system of matrix units for $M_n$.
Define a linear function 
$\ph : M_n \otimes e_1 A e_1 \rightarrow C^*(G, A, \alpha)$ 
by $\ph (e_{g,h} \otimes a ) = z_g a z_h^*.$  
One can then check in the usual way that $\ph$ is a homomorphism. 

It is also worth computing at this stage that for $g, h \in G$ and $ a \in e_1 A e_1$, 
we have $\| \ph (e_{g,h} \otimes a) -w_{g,1} a w_{h,1}^*\| \leq 2 \|a\| \delta_0$.
Let $D$ be the image of $\ph$, so that $\ph $ is clearly surjective as a map from $M_n \otimes e_1 A e_1$ to $D$. 
To check that $\ph $ is injective we first recall that 
$\ker (\ph) \cap (e_{g,h} \otimes e_1 A e_1)=e_{g,h} \otimes I$ 
where $I$ is an ideal of $e_1 A e_1$ which does not change 
as $g$ and $h$ vary.  But 
if $0=\ph (e_{g,h} \otimes a)=z_g a z_h^*$ for some 
$a \in e_1A e_1$, then multiplying on the left by $z_g^*$ 
and on the right by $z_h$ we see that $e_1 a e_1=a =0$, 
so $I=0$, thus $\ph$ is injective.

Now $\ph (e_{1,1}  \otimes a) = z_1 a z_1^* =e_1 a e_1 =a$ 
for any $a \in e_1 A e_1$.  Also, \\
$\ph (e_{g,g} \otimes 1)=z_g e_1 z_g^*=z_g z_g^* z_g z_g^*=e_g\in A$.  
These two conditions make up (1) of the conclusion.

To verify (2), let $a \in e_1Ae_1$ and estimate
\[ \|\ph (e_{g,g} \otimes a)- \alpha_g (a)\| \leq 2 \|a \| \delta_0 +\| u_g e_1 a e_1 u_g^* - \alpha_g (a) \| \leq \ep \|a\|.\]

For (4) we observe 
\[\sum_{g \in G} \ph (e_{g,g} \otimes 1) 
=\sum_{g\in G} z_g e_1 z_g^* =\sum_{g \in G} e_g =e.\]

Condition (5) holds essentially by construction since 
$1-e$ is Murray-von Neumann equivalent to a projection in 
$q_0 A q_0$, but $q_0 \in q A q$ and $q$ is equivalent 
to a projection in the hereditary subalgebra generated by $z$. 
 In total this gives $1-e$ is subequivalent to a projection 
in the hereditary subalgebra generated by $z$.

Now for condition (6), since $q_j \sim q_i$  we have 
$\tau (q_j) < \frac{1}{2N}$ for $0\leq j \leq 2N$ 
and for any $\tau \in T(A)$. In particular, since $1-e$ is subequivalent 
to $q_0$ we have $\tau (1-e) \leq \tau (q_0) < \frac{1}{2N}$.  
This implies $\frac{1}{2} < \tau (e)$.
Additionally $\tau (q_j) \leq \frac{1}{2N}$ implies 
$\tau (\sum_{j=1}^N q_j) < \frac{1}{2}$.
Combining these statements gives 
$\tau( \sum_{j=1}^N q_j) < \tau (e)$ for all $ \tau \in T(A)$.
So since order on projections over $A$ is determined by traces, 
$\sum_{j=1}^N q_j \precsim e$.  Let $h\in A$ be a projection satisfying
$\sum_{j=1}^N q_j \sim h \leq e$ and let $s$ be a partial isometry with 
$s^*s =\sum_{j=1}^N q_j$ and $s s^* =h$.  
Let $h_j=s q_j s^*$ for $j = 1, \dots , N$.  One checks that $h_1, \dots h_N$  
are mutually orthogonal projections summing to $h$.  Furthermore since 
$h_j \leq h \leq e$ we have $h_j \leq e$.  Furthermore, $h_j \sim g_j$ 
via the partial isometry $s g_j$.  So now we have 
$1-e \precsim g_j \sim h_j$.  
Let $f_j$ be a projection such that $1-e \sim f_j \leq h_j$.  
Since $f_j \leq h_j$, and the $h_j$ and $h_i$ are orthogonal for $1 \leq i,j \leq N$, we see that  
$f_1, \dots , f_N$ are mutually orthogonal.  
Finally $f_j \leq h_j \leq e$ in A and $eAe \subset eDe$, so $f_1, \dots , f_N$ are the projections we desired.

In order to show (3) we will use the following claim.

Claim: If $y=\sum_{g \in G} a_g u_g$ with $a_g \in A$ and $\|a_g\|  \leq M$, 
and if $[e_g,a_h]=0$ for all $g, h \in G$, 
then there are $d_1,\; d_2 \in D$ such that 
$\|e y - d_1\|$, $\|y e- d_2\|<8 n^2 M \delta_0$.

Proof of claim: We can write 
\[ey=\sum_{g \in G} \sum_{h \in G} e_g a_h u_h= 
\sum_{g \in G} \sum_{h \in G} (e_g a_h e_g) (e_g u_h)\]
since $e_g$ and $a_h$ commute.  
Now we make a norm estimate involving one of the factors in the third expression for $ey$ using the fact that $z_g$ is a partial isometry :
\begin{align*}
&\left\| \ph \left(e_{g,g} \otimes e_1 \alpha_{g^{-1}}\left(a_h\right) e_1\right) - e_g a_h e_g\right\| \\ 
&= \left\| e_g z_g \alpha_g^{-1}(a_h) z_g^* e_g - e_g a_h e_g\right\| \\ 
&\leq \left\|z_g \alpha_g^{-1}(a_h) z_g^* - z_g \alpha_g^{-1}(a_h) w_{g,1}^*\right\|
\\ & \hskip.5 in+ \left\|z_g \alpha_g^{-1}(a_h) w_{g,1}^* - w_{g,1}  \alpha_g^{-1}(a_h) w_{g,1}^*\right\|  
 + \left\|u_g e_1 \alpha_g^{-1}(a_h) e_1 u_g^*-e_g a_h e_g\right\| \\ 
& \leq 2M \delta_0 + 2M \delta.
\end{align*}

Now we make an estimate involving the other factor:
\begin{align*}
&\left\| \ph \left(e_{g,h^{-1}g} \otimes e_1\right) - e_g u_h\right\| \\
&\leq \left\| \ph \left(e_{g,h^{-1} g} \otimes e_1\right) - u_g e_1 u_{h^{-1} g}^*\right\| 
+ \left\|u_g e_1 u_{h^{-1} g}^*- e_g u_h\right\| \\
& \leq 2 \dt_0+ \dt.
\end{align*}

Let $d_0(g,h)= \ph (e_{g,g} \otimes e_1 \alpha_{g^{-1}}(a_h) e_1) \ph (e_{g,h^{-1}g} \otimes e_1)$.  
Then we have
\begin{align*}
 & \left\| d_0(g,h) - e_g a_h u_h \right\| \\ &
 \leq \left\| d_0(g,h) - \ph \left(e_{g,g} \otimes e_1 \alpha_{g^{-1}}(a_h) e_1\right) e_g u_h\ \right\| \\&
\alignInd + \left\|\ph \left(e_{g,g} \otimes e_1 \alpha_{g^{-1}}(a_h) e_1\right) e_g u_h -\left(e_g a_h e_g \right)\left(e_g u_h\right)\right\| \\ &
\leq \left\|\ph \left( e_{g,g} \otimes e_1 \alpha_{g^{-1}}(a_h) e_1\right)\right\| (2 \delta_0 +\delta) + 2M \delta_0 + 2M \delta \\ &
\leq 4M \delta_0 + 3M\delta.
\end{align*}

Now let $d_1= \sum_{g\in G} \sum_{h \in G} d_0(g,h)$.  Then
\begin{align*}
\| d_1 - ey \|&=\left \| \sum_{g\in G} \sum_{h \in G} d_0(g,h)- \sum_{g\in G} \sum_{h \in G} e_g a_h u_h \right\| \\
& \leq n^2 M (4 \delta_0+ 3\delta) \\&
< 8n^2 M \delta_0.
\end{align*}  

We now turn our attention to the construction of $d_2$.  
We can write 

$$ye=\sum_{g\in G} \sum_{h \in G} a_h \af_h (e_g) u_h.$$  

We note that 
\[\left\| \sum_{g\in G} \sum_{h \in G} a_h \af_h (e_g) u_h - \sum_{g\in G} \sum_{h \in G} a_h e_{hg} u_h \right\| <n^2 \dt . \]
 But  
\[ \sum_{g\in G} \sum_{h \in G} a_h e_{hg} u_h
 =\sum_{h\in G} \sum_{g \in G} e_{hg} a_h u_h
 =\sum_{h\in G} \sum_{k \in G} e_{k} a_h u_h\] 
 by making the change of variables, $k=hg$.  
 This last is of the same form as $ey$, so using the argument above 
 there is an element $d_2 \in D$ such that 
 \[\left\|\sum_{h\in G} \sum_{k \in G} e_{k} a_h u_h - d_2 \right\| \leq n^2 M (4 \delta_0+ 3M\delta).\]  
Thus 
\[\| y e - d_2\| \leq n^2 M (4 \delta_0+ 3M\delta) +n^2 \delta < 8n^2 M \delta_0.\]

We are now in a position to prove (3).  
Let $x \in F$ and choose $b_g \in S$ such that $x= \sum_{g \in G} b_g u_g$.
Define 
\[a_g = (1-e) b_g (1-e) + \sum_{h \in G} e_h b_g e_h . \] 
Now we get 
$$b_g -a_g = \sum_{h \in G} \left[(1-e) b_g e_h +e_h b_g (1-e)\right] 
+ \sum_{h\in G} \sum_{\stackrel{k \in G}{k\not= h}} e_k b_g e_h.$$

Because $b_g \in S$ we have,
$\left\| [b_g, e_h] \right\| < \delta $.  Then 
\begin{align*}
 \left\|b_g -a_g \right\| & 
\leq \sum_{h \in G} \left[ \left\|\left(1-e\right)\left[b_g, e_h\right] \right\| + \left\|\left(1-e\right)e_h b_g \right\| 
	+ \left\| \left[b_g, e_h\right]\left(1-e\right) \right\| + \left\|b_g e_h \left(1-e\right)\right\|  \right. \\ &
	\alignInd + \sum_{\stackrel{k \in G}{k\not= h}} \left\|e_k [b_g, e_h] \right\| 
	+\left. \sum_{\stackrel{k \in G}{k\not= h}} \left\|e_k e_h b_g \right\|  \right]\\ &
=(2n+n^2)\dt \\&
< 3 n^2 \delta .
\end{align*}

Set $y = \sum_{g \in G} a_g u_g$.  Then 

\[\|x-y \|  \leq \sum_{g\in G} \|(b_g-a_g) u_g \| 
\leq \sum_{g\in G} \|b_g -a_g\| \leq n (3n^2 \delta ) =3n^3 \delta .\]

One easily checks that $[a_g , e_k] =0 $ for all $ g, k \in G$.  Thus the claim applies to $y$ and provides $d_1 \in D$ such that 
$\|ey-d_1 \| <8n^2 M \delta_0$.  Therefore
\[\|ex -d_1 \| \leq \|ex-ey\| +\| ey -d_1 \| 
\leq 3n^3 \delta + 8n^2 M \delta_0
< \ep \]
 by the choice of $\delta$ and $\delta_0$.

Similarly, the claim provides $d_2 \in D$ such that 
$\|ye-d_2 \| < 8n^2 M \delta_0$ which implies that $d_2$ satisfies
$\|xe-d_2\| < \ep.$
\end{proof} 

Given objects satisfying part~(1) of the conclusion 
of Lemma~\ref{Basic1}, 
we can make a useful \hm\ into $C^* (G, A, \af)$ 
which should be thought of as 
a kind of twisted inclusion of $A.$ The following lemma 
is stated in terms of an arbitrary unital C*-algebra $B$, 
but we note it applies when $B=C^*(G,A, \af)$ and $\iota$ is the 
standard embedding.

\begin{lem}\label{EmbedA} 
Let $A$ be any simple unital \ca, let $B$ be a unital C*-algebra, and let 
$ \iota : A \rightarrow B$ be a unital injective homomorphism. 

Let $e, f \in A$ be \pj s, and let $n \in \N.$ 
Assume that there is an injective unital \hm\ %
$\ph \colon M_n \otimes f A f \to \io (e) B \io (e)$ 
such that, 
with $( e_{j, k} )$ being the standard system of matrix units for $M_n,$ 
we have $\ph (e_{1, 1} \otimes a) = \io (a)$ for all $a \in f A f.$ 
Then there is a corner $A_0 \S M_{n + 1} \otimes A$ which contains
\[ 
\left\{ 
\ts{ \left( \begin{array}{cc} a & 0 \\ 0 & b \end{array} \right) } 
\colon {\mbox{$a \in (1 - e) A (1 - e)$ and 
$b \in M_n \otimes f A f$}} \right\} 
\] 
as a unital subalgebra, and 
an injective unital \hm\ $\ps \colon A_0 \to B$ 
such that 
\[ 
\ps \left( \ts{ \begin{array}{cc} a & 0 \\ 0 & b \end{array} } \right) 
= \io (a) + \ph (b) 
\] 
for $a \in (1 - e) A (1 - e)$ and $b \in M_n \otimes f A f.$

Moreover, if $\af : G \rightarrow \Aut (A)$ is an action of a finite group on $A$, 
$B=C^*(G,A, \af)$, and $\iota $ is the standard inclusion, then for every 
$\af$-invariant tracial state $\ta$ on $A$ 
there is a tracial state $\sm$ on $C^* (G, A, \af)$ such that 
the extension ${\overline{\ta}}$ of $\ta$ to $M_{n + 1} \otimes A$ 
satisfies ${\overline{\ta}} |_{A_0} = \sm \circ \ps.$ 
\end{lem}

\begin{proof}
For the main part of the lemma we simply observe that the proof of Lemma 2.6 in \cite{OP1} works in this generality.  

For the part about traces, the proof of Lemma 2.6 in \cite{OP1} also works if one replaces the conditional expectation there with the conditional expectation $E(\sum_{g \in G} a_g u_g) =a_1$.
\end{proof} 


\section{TRACES AND ORDER ON PROJECTIONS IN CROSSED PRODUCTS}\label{Sec:OrdPj}

In this section, we prove that if 
$A$ is a simple unital C*-algebra with real rank zero such that 
the order on projections over $A$ is 
determined by traces, and if $\af: G \rightarrow \Aut (A)$ 
is an action of a finite group $G$ with  the \tRp, 
then the order on projections over $C^* (G, A, \af)$ is 
determined by traces. 

We begin with a comparison lemma for \pj s in 
crossed products by actions with the \tRp.

\begin{lem}\label{BasicOrder} 
Assume the hypotheses of Lemma~\ref{EmbedA} with $B=C^*(G,A, \af )$, 
and assume in addition that $A$ has real rank zero 
and that the order on projections over $A$ is 
determined by traces. 
Let $\ps \colon A_0 \to C^* (G, A, \af)$ 
be as in the conclusion of Lemma~\ref{EmbedA}. 
Suppose that $p, \, q \in \ps (A_0)$ are \pj s 
such that $\ta (p) < \ta (q)$ 
for all tracial states $\ta$ on $C^* (G, A, \af).$ 
Then there exists a \pj\ $r \in \ps (A_0)$ such that $r \leq q$ 
and $r$ is \mvnt\ to $p$ in $C^* (G, A, \af).$ 
\end{lem}

\begin{proof}
The same proof used for Lemma 3.1 of \cite{OP1} works by changing the group from $\Z$ to the finite group $G$.
\end{proof}

\begin{thm}\label{K0OrdDetByTraces} 
Let $A$ be an \idsuca\ with real rank zero, and 
suppose that the order on projections over $A$ is 
determined by traces. 
Let $\af: G \rightarrow \Aut (A)$ be an action of a finite group with the \tRp. 
Then the order on projections over $C^* (G, A, \af)$ is 
determined by traces. 
\end{thm}

\begin{proof}
We first observe that the hypotheses on $A$ imply that 
$A$ is finite, but $M_n(A)$ satisfies all the same hypotheses, 
so $A$ is in fact stably finite.

A special case of Lemma 3.9 of \cite{Ph30} implies that $\id_{M_n} \otimes \af $ as an action on $M_n \otimes A$ has the 
\tRp \ for any $n \in \N$.  Therefore, it suffices to consider projections in 
$C^*(G,A, \af )$.  

The rest of the proof can be done in the same way as the proof of Theorem 3.5 of \cite{OP1}.  The only changes are to use Lemma \ref{Basic1} instead of Lemma 2.5 of \cite{OP1} and Lemma \ref{EmbedA} instead of Lemma 2.6 of \cite{OP1}.
\end{proof} 


\section{REAL RANK OF CROSSED PRODUCTS}\label{Sec:RR}

\indent The main theorem of this section is \ref{chp3:thm1}.

\begin{thm}\label{chp3:thm1} 
Let $A$ be an \idsfsuca\ with real rank zero. 
Suppose that the order on projections over $A$ is 
determined by traces and $\af :G \rightarrow \Aut (A)$ 
is an action of a finite group with the \tRp. 
Then $C^* (G, A, \af)$ has real rank zero. 
\end{thm}

\begin{proof} 

Set $B = C^* (\Z, A, \af).$ 

As in the proof of Theorem~\ref{K0OrdDetByTraces}, 
the other hypotheses imply that $A$ is stably finite.

This proof is very similar to the proof of Theorem 4.5 of \cite{OP1} except simpler.  Therefore, instead of giving the details of the proof we will describe how the two proofs differ.  The norm estimates come out a bit different (the estimates can be improved), but we will not discuss that.

To justify the fact that the crossed product is simple, one should use Corollary 1.6 of \cite{Ph30}, instead of Corollary 1.14 of \cite{OP1}.  Similarly, in the proof of the integer case, Osaka and Phillips applied Lemma 2.5 of \cite{OP1} to obtain 
$e, \; p, $ and $f \in A$ projections, integers $n,m >0$, a unital subalgebra $D \subset eBe$ and an isomorphism $\varphi: D \rightarrow M_n \otimes fAf$.  
We apply Lemma \ref{Basic1} to obtain a \pj\ $e \in A \S C^* (G, A, \af),$ 
a unital subalgebra $D \S e C^* (G, A, \af) e,$ 
a \pj\ $f \in A,$ and an isomorphism 
$\ph \colon M_n \otimes f A f \to D.$

We use our $e$ to replace both $e$ and $p$ in the proof of Theorem 4.5 of \cite{OP1}.  Let $q=1-e$.  Since $a_0 e=e a_0$, there is no need to replace $p$ by a smaller projection approximately commuting with $a_0$.  This is the source of the improved estimates later on.
\end{proof}

\begin{cor}\label{Traces} 
Let $A$ be an \idsfsuca\ with real rank zero, 
and suppose that the order on projections over $A$ is 
determined by traces. 
Let $\af : G \rightarrow \Aut (A)$ be an action of a finite group with the \tRp. 
Then the restriction map is a bijection from the 
tracial states of $C^* (G, A, \af)$ to the $\af$-invariant 
tracial states of $A.$ 
\end{cor}

\begin{proof} 
This follows from Proposition~2.2 of~\cite{Ks3}, since $C^* (G, A, \af)$ has real rank zero by Theorem~\ref{chp3:thm1}.
\end{proof}


\section{STABLE RANK OF CROSSED PRODUCTS}\label{Sec:tsr}

\indent 
In this section, we prove that if 
$A$ is an infinite dimensional simple unital C*-algebra with real rank zero and 
stable rank one, such that 
the order on projections over $A$ is 
determined by traces, and if $\af : G \rightarrow \Aut (A)$ 
is an action of a finite group with  the \tRp, 
then $C^* (G, A, \af)$ has stable rank one.

\begin{lem}\label{LemmaForTSR} 
Let $A$ be an infinite dimensional simple unital C*-algebra 
with real rank zero and such that 
the order on projections over $A$ is determined by traces. 
Let $\af : G \rightarrow \Aut (A)$ be an action of a 
finite group with the \tRp. 
Let $q_1, \ldots, q_n \in C^* (G, A, \af)$ be nonzero \pj s, 
let $a_1, \ldots, a_m \in C^* (G, A, \af)$ be arbitrary, 
and let $\ep > 0.$ 
Then there exists a unital subalgebra $A_0 \S C^* (G, A, \af)$ 
which is stably isomorphic to $A,$ 
a \pj\ $p \in A_0,$ nonzero \pj s $r_1, \ldots, r_n \in p A_0 p,$ 
and elements $b_1, \ldots, b_m \in C^* (G, A, \af),$ 
such that: 
\begin{enumerate} 
\item 
$\| q_k r_k - r_k \| < \ep$ for $1 \leq k \leq n.$ 
\item
For $1 \leq k \leq n$ there is a \pj\ $g_k \in r_k A_0 r_k$ 
such that $1 - p \sim g_k$ in $C^* (G, A, \af).$ 
\item
$\| a_j - b_j \| < \ep$ for $1 \leq j \leq m.$ 
\item 
$p b_j p \in p A_0 p$ for $1 \leq j \leq m.$ 
\end{enumerate} 
\end{lem}

\begin{proof} 
The same argument used in the proof of Lemma 5.2 of \cite{OP1} works here as well.  As in the proof of 
\ref{chp3:thm1}, when we apply Lemma \ref{Basic1}, we set $p=e$ in order to use the proof of Lemma 5.2 of \cite{OP1}. 
\end{proof} 

\begin{thm}\label{chp4:thm1} 
Let $A$ be an infinite dimensional simple unital C*-algebra 
with real rank zero and stable rank one, and such that 
the order on projections over $A$ is determined by traces. 
Let $\af:G \rightarrow \Aut (A)$ be an action of a finite group with the \tRp. 
Then $C^* (G, A, \af)$ has stable rank one. 
\end{thm}

\begin{proof} 
The same proof used to show Theorem 5.3 of \cite{OP1} can be used to prove this theorem as well.  The only difference is instead of using Lemma 5.2 of \cite{OP1}, we must use Lemma \ref{LemmaForTSR}.
\end{proof} 


	
	\renewcommand{\bibname}{\textsc{REFERENCES}} 
\bibliographystyle{elsart-num-sort}	
	\addcontentsline{toc}{chapter}{\bibname}
	\bibliography{refs}					


\end{document}